\documentclass[11pt]{article}
\usepackage{amsmath, amssymb, amscd, amsthm, amsfonts}
\usepackage{graphicx}
\usepackage{algorithm,algorithmic}
\usepackage{hyperref}
\usepackage{tikz}
\usepackage{caption}
\usepackage{subcaption}
\usepackage{natbib}
\usepackage{bm}
\usepackage{multirow}
\usepackage{float}
\usepackage{titlesec}
\usepackage{authblk}

\titleformat{\paragraph}
{\normalfont\normalsize\bfseries}{\theparagraph}{1em}{}
\titlespacing*{\paragraph}
{0pt}{3.25ex plus 1ex minus .2ex}{1.5ex plus .2ex}

\title{Efficient proximal gradient algorithms for joint graphical lasso}
\author[1]{Jie Chen}
\author[1]{Ryosuke Shimmura}
\author[1]{Joe Suzuki}
\affil[1]{Graduate School of Engineering Science, Osaka University, Japan.}
\date{}

\newtheorem{theorem}{Theorem}
\newtheorem{lemma}[theorem]{Lemma}
\newtheorem{definition}{Definition}
\newtheorem{proposition}{Proposition}

\begin{document}

\maketitle
\begin{abstract}
We consider learning an undirected graphical model from sparse data. While several efficient algorithms have been proposed for graphical lasso (GL), the alternating direction method of multipliers (ADMM) is the main approach taken concerning for joint graphical lasso (JGL). We propose proximal gradient procedures with and without a backtracking option for the JGL. These procedures are first-order and relatively simple, and the subproblems are solved efficiently in closed form. We further show the boundedness for the solution of the JGL problem and the iterations in the algorithms.
The numerical results indicate that the proposed algorithms can achieve high accuracy and precision, and their efficiency is competitive with state-of-the-art algorithms. 

\end{abstract}

\section{Introduction}

Graphical models are widely used to describe the relationships among interacting objects \citep{lauritzen1996graphical}. 
Such models have been extensively used in various domains, such as bioinformatics, text mining, and social networks. 
The Graph provides a visual way to understand the joint distribution of an entire set of variables. 

In this paper, we consider learning Gaussian graphical models that are expressed by undirected graphs. 
It represents the relationship among continuous variables that follow a joint Gaussian distribution. 
In an undirected graph $\mathcal{G}=(V, E)$, edge set $E$ represents the conditional dependencies among the variables in vertex set $V$. 

Let $X_1,\ldots,X_p$ ($p\geq 1$) be Gaussian variables with covariance matrix $\bm{\Sigma}\in {\mathbb R}^{p\times p}$, and $\bm{\Theta}:=\bm{\Sigma}^{-1}$ be, if it exists, the precision matrix.
We draw the edges so that the variables $X_i, X_j$ are conditionally independent given the other variables if and only if the $(i,j)$-th element $\theta_{i, j}$ in $\bm{\Theta}$ is 0:
\[
\{i,j \}\not\in E\ \ \Longleftrightarrow\  \theta_{i,j}=0 \ \  \Longleftrightarrow  X_i \perp \!\!\! \perp X_j |X_{V\backslash \{i,j\} }\ ,
\]
where each edge is expressed as a set of two elements in $\{1,\ldots, p\}$.
%It can be shown to correspond to the underlying graph structure. 
In this sense, constructing a Gaussian graphical model is equivalent to estimating a precision matrix. 

Suppose that we estimate the undirected graph from data consisting of $n$ tuples of $p$ variables and that dimension $p$ is much higher than sample size $n$.
For example, if we have expression data of $p=20,000$ genes for $n=100$ case/control patients, how can we construct a gene regulatory network structure from the data?
However, it is almost impossible to estimate the locations of the nonzero elements in $\bm{\Theta}$ by obtaining the inverse of the estimate $S\in {\mathbb R}^{p\times p}$ of the covariance matrix $\bm{\Theta}$. 
In fact, if $p>n$, then no inverse $S^{-1}$ exists because the rank of $\bm{S}\in {\mathbb R}^{p\times p}$ is, at most, $n$.

To address the situation, two directions are suggested:
\begin{enumerate}
\item Sequentially find the variables on which each variable depends via regression so that the quasi-likelihood is maximized \citep{meinshausen2006high}
\item Find the locations in $\bm{\Theta}$, the values of which are zeros, so that the $\ell_1$ regularized log-likelihood is maximized \citep{yuan2007model,friedman2008sparse,banerjee2008model,rothman2008sparse} %(graphical Lasso, GL).
\end{enumerate}
We follow the second approach because we assume Gaussian variables, also known as graphical lasso (GL). The 
$\ell_1$ regularized log-likelihood is defined by 
\begin{align} \label{GL}
\max_{\bm{\Theta}}\log\det \bm{\Theta} - \text{trace}S \bm{\Theta} -\lambda || \bm{\Theta} ||_1
\end{align}
where tuning parameter $\lambda$ controls the amount of sparsity, and  
$||\bm{\Theta}||_1$ denotes the sum of the absolute value of the off-diagonal elements in $\bm{\Theta}$. Several optimization techniques \citep{banerjee2006convex,mazumder2012graphical,guillot2012iterative,d2008first,friedman2008sparse,hsieh2014quic} have been studied for the optimization problem of \eqref{GL}.

In particular, we consider a generalized version of the abovementioned GL.
For example, suppose that the gene regulatory networks of thirty case and seventy control patients are different.
One might construct a gene regulatory network separately for each of the two categories.
However, estimating each on its own does not provide an advantage if a common structure is shared. 
Instead, we use one hundred samples to construct two networks simultaneously.
Intuitively speaking, using both types of data improves the reliability of the estimation by increasing the sample size
for the genes that show similar values between case and control patients,  
while using only one type of data leads to a more accurate estimate for genes that show significantly different values.
\citet{danaher2014joint} proposed a joint graphical lasso (JGL) model by including an additional convex penalty (grouped lasso penalties) 
to the graphical lasso objective function for $K$ classes. For example, $K$ is two for the case/control in the example. 
Although there are several approaches to the problem, such as 
\cite{honorio2010multi}, \cite{guo2011joint}, \cite{zhang2012learning}, \cite{hara2013learning},
the JGL is considered the most promising.

The main topic of this paper is to improve efficiency in terms of solving the JGL problem.
For the GL, a relatively efficient solving procedures exists.
If we differentiate the $\ell_1$ regularized log-likelihood by $\bm{\Theta}$, then we have an equation 
to solve \citep{friedman2008sparse}.
Moreover, several improvements have been considered for the GL, such as proximal Newton \citep{hsieh2014quic} and proximal gradient \citep{guillot2012iterative} procedures.
However, for the JGL, even if we derive such an equation, we have no efficient way to handle it.%

Instead,  the alternating direction method of multipliers (ADMM) \citep{glowinski1975approximation},
which is a procedure for solving convex optimization problems for general purposes,
 has been the main approach taken \citep{danaher2014joint,tang2015exact,hallac2017network,gibberd2017regularized}. However, ADMM  does not scale well concerning feature dimension $p$ and number of classes $K$. It usually takes time for convergence to obtain high accuracy \citep{boyd2011distributed}.
 
\begin{table}
\caption{Efficient JGL Procedures \label{tab01}}
\begin{center}
\begin{tabular}{l|l|l|l}
\hline
GL/JGL&Original&Proximal Newton&Proximal Gradient\\
\hline
GL& \cite{friedman2008sparse}& QUIC \citep{hsieh2014quic}  & \cite{guo2011joint}\\
\hline
JGL& \cite{danaher2014joint}& \cite{yang2015fused} (for fused penalty)&{\color{red}Current Paper}\\
\hline
\end{tabular}
\end{center}
\end{table}

For the efficient procedures of the JGL problem, \citet{yang2015fused} proposed a method based on the proximal Newton procedure only when the penalty term is expressed by fused lasso (FMGL). 
The existing method requires expensive computations for the Hessian matrix and Newton directions, which means that it would be expensive for high-dimensional problems. 

In this paper, we propose efficient proximal-gradient-based algorithms to solve the JGL problem
by extending the procedure in \cite{guillot2012iterative} and modifying the step-size selection strategy proposed in \cite{tran2015composite}.  
Moreover, we provide the theoretical analysis of both methods for the JGL problem.

We show that the proposed methods are faster than ADMM and FMGL for any class $K$ and large-scale $p$.
Moreover, in our proximal gradient method for the JGL problem, 
the proximal operator in each iteration is quite simple, which eases the implementation process and requires very little computation and memory at each step. 
Simulation experiments are used to justify our proposed methods over the existing ones.

Our main contributions are as follows:
\begin{itemize}
	\item We propose efficient algorithms based on the proximal gradient method to solve the JGL problem. The algorithms are first-order and quite simple, and the subproblems can be solved efficiently with a closed-form solution.  The numerical results indicate that the methods can achieve high accuracy and precision, and the computational time is competitive with that of state-of-art algorithms.%Our algorithms are computationally efficient with competitive performance over the state of the arts.
	\item We provide the boundedness for the solution of the JGL problem and the iterations in algorithms. Then, the domain for the precision matrix in each iteration in our algorithms is constrained, guaranteeing the iterations inside the compact domain. 
\end{itemize} 
%Moreover, there are theoretical guarantees for the proposed algorithms by the help of \cite{guillot2012iterative} and \cite{tran2015composite}.
Table \ref{tab01} summarizes the relationship between the proposed and existing methods.

The remaining parts of this paper are as follows. In Section 2, we first provide the background of our proposed method and introduce the joint graphical lasso problem. In Section 3, we illustrate the detailed content of the proposed algorithms and provide some theoretical analysis. In Section 4, we report some numerical results of the proposed approaches, including comparisons with efficient methods and performance evaluations. Finally, we draw some conclusions in Section 5.

\section{Preliminaries}

This section first reviews the graphical lasso (GL) problem and the G-ISTA algorithm \citep{guillot2012iterative} to solve the GL problem. 
Then, we introduce the step-size selection strategy that we employed and the details of the joint graphical lasso problem.
\subsection{Graphical lasso}
Let $\bm{x}_1,\ldots,\bm{x}_n \in \mathbb{R}^p$ be $n\geq 1$ observations of dimension $p\geq 1$ that follow the Gaussian distribution with mean $\bm{\mu}\in {\mathbb R}^{p}$ and covariance matrix $\bm{\Sigma} \in {\mathbb R}^{p\times p}$, where without loss of generality, we assume $\bm{\mu} = \bm{0}$. Let $\bm{\Theta} = \bm{\Sigma}^{-1}$, and the empirical covariance matrix $\bm{S}:=\frac{1}{n}\sum_{i=1}^n \bm{x}_i^T\bm{x}_i$.
Given penalty parameter $\lambda>0$, 
the graphical lasso (GL) is the procedure to find the positive definite $\bm{\Theta}\in {\mathbb R}^{p\times p}$:
\begin{equation}\label{eq1}
\underset{\bm{\Theta}}{\text{minimize}}  -\log\det\bm{\Theta}+{\rm trace}(\bm{S\Theta})+\lambda\|\bm{\Theta}\|_1\ ,
\end{equation}
where $||\bm{\Theta}||_1= \sum_{j\neq k} |\theta_{j,k}|. $ If we regard $V:=\{ 1,\dots, p\}$ as a vertex set, then we can construct an undirected graph with edge set $\{\{j,k\}|\theta_{j,k} \neq 0  \}$, where set $\{j,k\}$ denotes an undirected edge that connects the nodes $j,k \in V$.

If we take the subgradient of \eqref{eq1}, then we find that the optimal solution $\bm{\Theta}_*$ satisfies the optimality condition:
\begin{equation}\label{eq2}
\bm{\Theta}_*^{-1}-\bm{S}-\lambda\bm{\Phi}=0\ ,
\end{equation}
where $\bm{\Phi}=(\Phi_{j,k})$ is
$$\Phi_{j,k}=
\left\{
\begin{array}{ll}
1,&\theta_{j,k}>0\\
{[-1,1]},&\theta_{j,k}=0\\
-1,&\theta_{j,k}<0
\end{array}
\right.\ .
$$

\subsection{ISTA for graphical lasso} \label{istagl}

In this section, we introduce the procedure for solving the GL problem \eqref{eq1} by the iterative shrinkage-thresholding algorithm (ISTA) proposed by \citet{guillot2012iterative}, which is a proximal gradient method usually employed in dealing with nondifferentiable optimization problems. 

Specifically, the general ISTA solves the following composite optimization problems:
\begin{equation} \label{eq2.2}
	\underset{\bm{x}}{\text{min }} F(\bm{x}):= f(\bm{x}) +g(\bm{x})
\end{equation}
where $f$ and $g$ are convex functions, with $f$ continuously differentiable and $g$ possibly nonsmooth.

 For the GL problem \eqref{eq1}, denote $f, g:\mathbb{R}^{p\times p} \rightarrow \mathbb{R}$ as
$$f(\bm{\Theta}):=-\log\det\bm{\Theta}+{\rm trace}(\bm{S\Theta})\ ,$$
and 
$$g(\bm{\Theta}):=\lambda\|\bm{\Theta}\|_1\ .$$
If we define the quadratic approximation $Q_{\eta}:\mathbb{R}^{p\times p} \times \mathbb{R}^{p\times p} \rightarrow \mathbb{R}$ w.r.t. $f(\bm{\Theta})$:
\begin{align} \label{eqbc}
	Q_\eta(\bm{\Theta}',\bm{\Theta}) := f(\bm{\Theta})+<\bm{\Theta}'-\bm{\Theta},\nabla f(\bm{\Theta})>+\frac{1}{2\eta} ||\bm{\Theta}'-\bm{\Theta}||^2_F\ ,
\end{align}
then we can describe the ISTA as a procedure that iterates
\begin{align} \label{eq87}
	\bm{\Theta}_{t+1} &= \arg\min_{\bm{\Theta}} \{Q_{\eta_t} (\bm{\Theta},\bm{\Theta}_t) +g(\bm{\Theta}) \} \\
	&=  {\rm prox}_{\eta_t g}(\bm{\Theta}_t - \eta_t \nabla f ( \bm{\Theta}_t) )
\end{align}
given initial value $\bm{\Theta}_0$, where the  value of step size $\eta_t$ may change at each iteration $t=1,2,\dots,$ and use the notation of the proximity operator:
$${\rm prox}_h(\bm z):={\rm argmin}_{\theta}\{\frac{1}{2}\|\bm z-\bm \theta\|_2^2+h(\bm \theta)\} . $$
Note that the proximal operator of function $g=\lambda ||\bm \Theta||_1 $ is the soft-thresholding operator: the absolute value $|\theta_{i,j}|$ of each off-diagonal element $\theta_{i,j}$ with $i\not=j$ becomes either $|\theta_{i,j}|-\lambda$ or zero (if $|\theta_{i,j}|<\lambda$). 
We use the following function in Sections 2.3 and 3:
\begin{equation}\label{eq65}
{\cal S}_\lambda (\bm{\Theta}) =\text{sgn}(\theta_{i,j}) (|\theta_{i,j}-\lambda|)_{+}
\end{equation}
where $(x)_+ := max(x,0)$.

\begin{definition} %cite
A differentiable function $f: \mathbb{R}^n\rightarrow \mathbb{R}$ is said to have a Lipschitz-continuous gradient if there exists $L>0$ (Lipschitz constant) such that
\begin{equation} \label{lip}
||\nabla f(\bm{x})-\nabla f(\bm{y})|| \leq L ||\bm{x}-\bm{y}||,\forall \bm{x},\bm{y} \in \mathbb{R}^n	
\end{equation}
\end{definition}

It is known that if we choose $\eta_t=1/L$ for each step in the 
ISTA that minimizes $F(\cdot)$, then the convergence rate is, at most, 
\begin{equation}\label{eq77}
F(\bm{\Theta}_t)-F(\bm{\Theta}^*)=O(1/t)
\end{equation}
\citep{beck2009fast}. However, for the GL problem \eqref{eq1}, we know neither the exact value of the Lipschitz constant $L$ nor any nontrivial upper bound. \citet{guillot2012iterative} implement a backtracking line search option in step 1 of Algorithm \ref{alg:algorithm0} below to handle this issue.

The backtracking line search means that we compute the $\eta_t$ value for each time $t=1,2,\ldots$ by repeatedly multiplying $\eta_t$ by  $0<c<1$ until $\bm{\Theta}\succ 0$ ($\Theta$ is positive definite) and 
\begin{equation}\label{eqbck}
f(\bm{\Theta}_{t+1})\leq Q_{\eta_t}(\bm{\Theta}_{t+1},\bm{\Theta}_t)
\end{equation}
for the $\bm{\Theta}_{t+1}$ in  (\ref{eq87}),
which means that $f(\bm{\Theta}_{t})$ is nonincreasing:
$$f(\bm{\Theta}_{t+1})\leq Q_{\eta_t}(\bm{\Theta}_{t+1},\bm{\Theta}_t)\leq
Q_{\eta_t}(\bm{\Theta}_{t},\bm{\Theta}_t)=f(\bm{\Theta}_{t})\ .
$$
Additionally, (\ref{eqbck}) is a sufficient condition for (\ref{eq77}), which was derived in \cite{beck2009fast} (see the relation between Lemma 2.3 and Theorem 3.1).

Moreover, the smaller $\eta_t$ is, the more likely it satisfies condition \eqref{eqbck}. $\bm{\Theta}_{t+1}$ is continuous w.r.t. $\eta_t$ when $\bm{\Theta}_t$ is fixed, which means that $\bm{\Theta}_{t+1}$ trivially satisfies the 
condition for $\eta_t=0$. 
%Hence, the upper bound of L (i.e. $\frac{1}{\eta_t} \geq  L>0$) also make sense.

The whole procedure is given in the following algorithm, and we continue describing the third and fourth steps.

\begin{algorithm}[H]
\caption{G-ISTA for problem \eqref{eq1} }
\label{alg:algorithm0}
\textbf{Input}: S, tolerance $\epsilon>0$, backtracking constant $0<c<1$, initial value $\eta_0$, $\bm{\Theta}_0$.  \\
\textbf{While} $\Delta>\epsilon$ \textbf{do} 
\begin{algorithmic}[1] %[1] enables line numbers
\STATE backtracking line search: \\
While
\begin{equation} \label{condition}
	\bm{\Theta}_{t+1}\succ  0  \quad \text{and} \quad f(\bm{\Theta}_{t+1}) \leq Q_{\eta_{t}} (\bm{\Theta}_{t+1},\bm{\Theta}_t)
\end{equation}
do not hold for $\bm{\Theta}_{t+1}:= {\cal S}_{\eta_{t} g} (\bm{\Theta}_{t}- \eta_{t} \nabla f(\bm{\Theta}_{t}))$ 
multiply $\eta_t$ by c.
\STATE Update iterate: $\bm{\Theta}_{t+1}\leftarrow {\cal S}_{\eta_{t} g} (\bm{\Theta}_{t}- \eta_{t} \nabla f(\bm{\Theta}_{t}))$
\STATE Set next initial step $\eta_{t+1}$ by the Barzilai-Borwein gradient.
\STATE Compute duality gap: $\Delta$
\end{algorithmic}

\textbf{Output}: $\epsilon$-optimal solution to problem (1), $\bm{\Theta}^*=\bm{\Theta}_{t+1 }$
\end{algorithm}

The following properties provide a theoretical foundation that guarantees efficient convergence and suggest an efficient strategy for step size in the backtracking procedure and step 3.
\begin{lemma} (\cite{guillot2012iterative}, lemma 3) \label{lm1}
Let $\{\bm{\Theta}_t \}_{t=0,1,\cdots}$ be the sequence generated by Algorithm \ref{alg:algorithm0}, and $\bm{\Theta}_*$ be the optimal solution of the problem \eqref{eq1}; moreover, let
\[
a:=\min\{\lambda_{\min}(\bm{\Theta}_*), \lambda_{\min}(\bm{\Theta}_t) \}
\]
and
\[
b:=\max\{\lambda_{\max}(\bm{\Theta}_*),  \lambda_{\max}(\bm{\Theta}_t) \}\\ , 
\]
where $\lambda_{\max}(\bm{\Theta}_*)$ and $\lambda_{\max}(\bm{\Theta}_*)$ are the maximum and minimum among the $p$ eigenvalues of $\bm{\Theta}_*$,
and $\lambda_{\max}(\bm{\Theta}_t)$ and $\lambda_{\max}(\bm{\Theta}_t)$ are the maximum and minimum among all the eigenvalues of $\bm{\Theta}_t$, $t=1,2,\ldots$, respectively. Then, we have 
\[
\|\bm{\Theta}_{t+1}-\bm{\Theta}^*\|_F  \leq \gamma \|\bm{\Theta}_t-\bm{\Theta}_{*}\|_F
\]
with the convergence rate
\[
\gamma:=\max\{ |1- \frac{\eta_t}{b^2}|,|1- \frac{\eta_t}{a^2}|\}\ \ .
\]
\end{lemma}
Lemma \ref{lm1}  implies that to obtain $\gamma<1$, we require
\begin{equation}\label{eq81}
0<\eta_t<a^2\ 
\end{equation}
If the step size is in the range of (\ref{eq81}), then the sequence $\{\bm{\Theta}_t\}_{t=0,1,\dots}$ is bounded as follows:
\begin{lemma}(\cite{hsieh2014quic}, Lemma2,  \cite{guillot2012iterative}, Theorem 2) \label{lm2}
	If the step size is in the range of (\ref{eq81}), then all the iterations $\{\bm{\Theta}_t\}_{t=0,1,\cdots,}$ belong to level set U defined as follows:
	\[
	U=\{ \bm{\Theta}_t | F(\bm{\Theta}_t) \leq F(\bm{\Theta}_0) \}
	\]
	 and  we have 
	\[
	m\leq \|\bm{\Theta}_t\| \leq M
	\]
	for the constant $m:=e^{-F(\bm{\Theta}_0)}M^{-(p-1)} $, $M:= ||\bm{\Theta}^*||_2+||\bm{\Theta}_0-\bm{\Theta}^*||_F$.
\end{lemma}

For each iteration of Algorithm \ref{alg:algorithm0}, the initial step of $\eta_t$ is needed. \citet{guillot2012iterative} considered step size $\eta := \lambda_{min}(\bm{\Theta}_t)^2$ using lemma \ref{lm1} at first. However, it is quite conservative in practice. Hence, they used the Barzilai-Borwein method \citep{barzilai1988two} for implementation and regarded the step size $\eta$ as safe. When the backtracking iterations in Step 1 exceed the given maximum number to satisfy condition \eqref{condition}, we can use the safe step for the subsequent calculations.

For Step 4, they considered the duality gap as a stopping criterion. The duality gap is defined as the difference between the primal and dual problems,  which implies optimality when close to zero. To illustrate the duality gap in Step 4, we first introduce the dual problem of the GL problem \eqref{eq1}.
\begin{align*}
&\min_{\bm{U}\in \mathbb{R}^{p\times p}} -\text{logdet}(\bm{S}+\bm{U})-p \\
&\text{subject to } \quad ||\bm{U}||_{\infty} \leq \lambda 	
\end{align*}
where $\bm{U}$ are the dual variables, and the primal and dual variables are connected by $\bm{\Theta}=(\bm{S}+\bm{U})^{-1}$.
Hence, the duality gap of the GL problem in Step 4 is
\[
    \Delta:= -\text{logdet}(\bm{S}+\bm{U}_{t+1})-p-\text{logdet} \bm{\Theta}_{t+1} 
    +<\bm{S},\bm{\Theta}_{t+1}> +\lambda ||\bm{\Theta}_{t+1}||_1\
\]

\subsection{Composite self-concordant minimization} \label{cscm}
\citet{tran2015composite} considered a composite version of self-concordant minimization, which was initially proposed by \citet{nesterov1994interior}.
\cite{tran2015composite} proposed a way to efficiently calculate the step size for the proximal gradient method for the GL problem. In particular, they proved that
$$f(\bm{\Theta}):=-\log\det\bm{\Theta}+{\rm trace}(S\bm{\Theta})\ $$ in \eqref{eq1} is self-concordant, and considers the minimization
\[
F^* := \min_{x} \{ F(x) := f(x) +g (x) \}
\] 
when $f$ is convex, and self-concordant and when $g$ is convex and nonsmooth. As for Algorithm \ref{alg:algorithm0}, 
without using the backtracking line search,
we can compute direction $d_t$:
\[
\bm{d}_t:={\cal S}_{L_{t}^{-1} g} (\bm{\Theta}_{t}- \eta_{t} \nabla f(\bm{\Theta}_{t}))- \bm{\Theta}_t\ ,
\]
where ${\cal S}_a$ with $a>0$ defined by (\ref{eq65}). Then, the step size can be determined by direction $d_t$. Let $L_t:=\eta_t^{-1},  \beta_t:=L_t||\bm{d}_t||_F^2 ,\lambda_t:=||\bm{\Theta}^{-1}_t \bm{d}_t||_F $, and the step size is 
\[
\alpha_t := \frac{\beta_t}{\lambda_t(\lambda_t +\beta_t)}\ ,
\]
after that update $\bm{\Theta}_{t+1}:= \bm{\Theta}_t +\alpha_t \bm{d}_t$ in the iterations.

In addition, \cite{tran2015composite} claimed that this strategy ensures a descent direction in the proximal gradient scheme and guarantees convergence.

\subsection{Joint Graphical Lasso}

Let $n\geq 1, p\geq 1, K\geq 2$, and $(\bm {x}_1,y_1),\ldots,(\bm{x}_n,y_n)\in {\mathbb R}^p\times \{1,\ldots,K\}$,
where each $x_i$ is a row vector. Let $n_k$ be the number of occurrences in $y_1,\ldots, y_n$ such that $y_i=k$, so that $\sum_{k=1}^Kn_k=n$.

For each $k=1,\ldots,K$, we define the empirical covariance matrix $\bm{S}^{(k)}\in {\mathbb R}^{p\times p}$ of the data $\bm{x}_i$ as follows:
$$\bm{S}^{(k)}:=\frac{1}{n_k}\sum_{i: y_i=k}\bm{x}_i^T \bm{x}_i\ $$
Given the penalty parameters $\lambda_1>0$ and $\lambda_2 >0 $, the joint graphical lasso (JGL) is the procedure to find the positive definite matrix $\bm{\Theta}^{(k)}\in {\mathbb R}^{p\times p}$ for $k=1,\ldots,K$, such that
\begin{equation}\label{eq6}
    \min_{\bm{\Theta}} -\sum_{k=1}^Kn_k\{\log\det\bm{\Theta}^{(k)}-{\rm trace}(S^{(k)}\bm{\Theta}^{(k)})\} \\
+\lambda_1 \sum_{k=1}^{K} \sum_{i \neq j} |\theta_{i,j}^{(k)}| + P(\bm{\Theta})\ ,
\end{equation}
where $P(\bm{\Theta})$ penalizes $\bm{\Theta}:=[\bm{\Theta}_1,\ldots,\bm{\Theta}_K]^T$, and the choices are diverse.
For example, \cite{danaher2014joint} suggested the fused and group lasso penalties: 
$$P_F(\bm{\Theta}):=\lambda_2 \sum_{k\not=l}\sum_{i,j}|\theta_{k,i,j}-\theta_{l,i,j}|$$
and 
$$P_G(\bm{\Theta}):=\lambda_2\sum_{i\not=j}\left\{\sum_{k=1}^K \theta_{k,i,j}^2\right\}^{1/2}\ ,$$
where $\theta_{k,i,j}$ is the $(i,j)$-th element of $\bm{\Theta}^{(k)}\in {\mathbb R}^{p\times p}$ for $k=1,\dots,K$.

\begin{lemma}(\cite{yang2015fused}) \label{uniqueness}
 Under the assumption that $\text{diag}(\bm{S}^{(k)})>0,k=1,\dots,K$, problem \eqref{eq6} has a unique optimal solution.
 \end{lemma}
 
Lemma \ref{uniqueness} shows the uniqueness of problem \eqref{eq6}. Unfortunately, there is no equation like (\ref{eq2}) for the JGL to find the optimum $\bm{\Theta}_*$. \cite{danaher2014joint} considered the ADMM to solve the JGL problem. However, ADMM is quite time-consuming for large-scale problems.

\section{Proposed Method}

In this section, we propose two efficient algorithms for solving the JGL problem. One is an extended ISTA based on the G-ISTA, and we modify the other by the step-size selection strategy illustrated in Section \ref{cscm}. 

\subsection{ISTA for JGL problem}
To neatly describe the JGL problem, we define $f,g:\mathbb{R}^{K\times p \times p} \rightarrow \mathbb{R}$ by
\begin{align} 
	f(\bm{\Theta})& := -\sum_{k=1}n_k\{\log\det\bm{\Theta}^{(k)}-{\rm trace}(S^{(k)}\bm{\Theta}^{(k)})\} \label{ftheta} \\
	g(\bm{\Theta})&: = \lambda_1 \sum_{k=1}^K \|\bm{\Theta}^{(k)}\|_1 + P(\bm{\Theta}) \label{gtheta}
\end{align}
Then, the problem \eqref{eq6} becomes the following:
\[
\min_{\bm{\Theta}} F(\bm{\Theta}):=f(\bm{\Theta})+g(\bm{\Theta})
\]
$f(\bm{\Theta})$ is convex and continuously differentiable, and $g(\bm{\Theta})$ is convex and nonsmooth. Therefore, the ISTA is available for solving the JGL problem \eqref{eq6}. 

The main differences between the G-ISTA and the proposed method are that the latter needs to consider $K$ categories of graphical models simultaneously in the JGL problem \eqref{eq6}. There are two combined penalties in $g(\Theta)$, which complicat the proximal operator in the ISTA procedure. Consequently, the operator for the proposed method is not a simple soft thresholding operator, as is that for the G-ISTA method.

If we define quadratic approximation $Q_{\eta_t}:\mathbb{R}^{K\times p\times p} \rightarrow \mathbb{R}$ of $f(\Theta)$:
\begin{align*}
Q_{\eta_t}(\bm{\Theta},\bm{\Theta}_t):=-\sum_{k=1}^Kn_k\{\log\det\bm{\Theta}^{(k)}-{\rm trace}(\bm{S}^{(k)}\bm{\Theta}^{(k)})\} \\
+\sum_{k=1}^K<\bm{\Theta}^{(k)}-\bm{\Theta}_t^{(k)},\nabla f(\bm{\Theta}_t^{(k)})> +\frac{1}{2\eta_t} \sum_{k=1}^K||\bm{\Theta}^{(k)}-\bm{\Theta}_t^{(k)}||_F^2\ ,    
\end{align*}

then the update iteration becomes the following:
\begin{equation*}
\bm{\Theta}_{t+1} =  \underset{\bm{\Theta}}{\text{argmin }} \left\{ Q_{\eta_t}(\bm{\Theta},\bm{\Theta}_t)+g(\bm{\Theta}) \right\}=  \text{prox}_{\eta_{t} g}
 (\bm{\Theta}_t - \eta_t \nabla f(\bm{\Theta}_t))    
\end{equation*}

Nevertheless, the Lipschitz gradient constant of $f(\Theta)$ is unknown over the whole domain in the JGL problem. Therefore, our approach also needs a backtracking line search to calculate step size $\eta_t$. We show the details in Algorithm 2.

\begin{algorithm}[H]
\caption{ISTA for problem \eqref{eq6} }
\label{alg:algorithm1}
\textbf{Input}: S, tolerance $\epsilon>0$, backtracking constant $0<c<1$, initial value $\eta_0$, $\bm{\Theta}_0$.  \\
\textbf{For} $t=0,1,\cdots,  $ (until convergence) \textbf{do} 
\begin{algorithmic}[1] %[1] enables line numbers
\STATE backtracking line search: 
While
\begin{align} \label{condi}
	f(\bm{\Theta}_{t+1}) \leq Q_{\eta_{t}} (\bm{\Theta}_{t+1},\bm{\Theta}_t)  \; \text{and} \; \bm{\Theta}_{t+1}^{(k)}\succ 0 \text{ for } k = 1,\cdots,K
\end{align}
do not hold for $\bm{\Theta}_{t+1}:= \text{prox}_{\eta_{t} g} (\bm{\Theta}_{t}- \eta_{t} \nabla f(\bm{\Theta}_{t}))$ 
multiply $\eta_t$ by c.
\STATE Update iterate: $\bm{\Theta}_{t+1}\leftarrow \text{prox}_{\eta_{t} g} (\bm{\Theta}_{t}- \eta_{t} \nabla f(\bm{\Theta}_{t}))$
\STATE Set next initial step $\eta_{t+1}$
\end{algorithmic}

\textbf{Output}: $\epsilon$-optimal solution to problem (1), $\Theta^*=\Theta_{t+1 }$
\end{algorithm}

After obtaining step size $\eta_t$, we discuss the detailed calculation of the proximal operator for fused lasso and group lasso penalties. In the following, we will see that the subproblems of the proximal operator for the two types of penalties are equivalent to the fused lasso and group lasso problems, respectively, and can be solved easily by any associated efficient procedure.

\subsubsection*{A. Fused lasso penalty $P_F$}
By the definition of the proximal operator in the update step, we have
\begin{equation}
 \bm{\Theta}_{t+1}  =	\arg\min_{\Theta} \left\{ \frac{1}{2}||\bm{\Theta}-\bm{\Theta}_{t}+\eta_{t} \nabla f(\bm{\Theta}_{t}) ||^2   + \eta_t \cdot \lambda_1 \sum_{k=1}^K \sum_{i\neq j}|\theta_{k,i,j}| + \eta_t \cdot \lambda_2 \sum_{k\not=l}\sum_{i,j}|\theta_{k,i,j}-\theta_{l,i,j}| \right\} \label{proxeq}
\end{equation}

Problem \eqref{proxeq} is separable with respect to the elements $\theta _{k,i,j}$ in $\Theta^{(k)} \in \mathbb{R}^{p\times p}$; hence, the proximal operator can be computed in componentwise operations:

Let $A=\bm{\Theta}_{t}-\eta_{t} \nabla f(\bm{\Theta}_{t}) $; then, problem \eqref{proxeq} reduces to the following subproblem for $i=1,\cdots,p$, $j=1,\cdots,p$ : 
\begin{equation}
  \underset{\theta_{1,i,j},\cdots,\theta_{K,i,j}}{\text{argmin }} \{ \frac{1}{2} \sum_{k=1}^K(\theta_{k,i,j}-a_{k,i,j} )^2 + \eta_t \cdot \lambda_1 \cdot 1_{i\neq j} \sum_{k=1}^K |\theta_{k,i,j}|+\eta_t \cdot \lambda_2 \sum_{k\not=l}|\theta_{k,i,j}-\theta_{l,i,j}| \} \label{flsa} 	
\end{equation}
where $1_{i\neq j}$ is an indicator function, the value of which is 1 only when $i\neq j$. Then, this is known as the fused lasso signal approximator (FLSA) \citep{friedman2007pathwise}. Several rather efficient algorithms can be used to solve this problem \citep{hoefling2010path,tibshirani2011solution,johnson2013dynamic}. 

In particular, for illustration, let $\alpha_1 := \eta_t \cdot \lambda_1 \cdot 1_{i\neq j}$ and $\alpha_2 := \eta_t \cdot \lambda_2$. 
When $i=j$, $\alpha_1 = 0, \alpha_2>0$,  subproblem \eqref{flsa} solves the FLSA problem with $K-$dimensional variable $\theta_{k,i,j}$.

When $i\neq j$, $\alpha_1 >0$, the solution to \eqref{flsa} can be obtained through soft thresholding based on the solution when $\alpha_1 =0$ by the following lemma.

  \begin{lemma} (\cite{friedman2007pathwise})
 Assume that the solutions to $\alpha_1 = 0 $ and $\alpha_2$ are known and denoted by $\theta_k (0,\alpha_2)$.
 	Then, the solution for $\alpha_1 >0$ in the fused lasso problem is
 	\[
 	\theta_k (\alpha_1,\alpha_2) = {\cal S}_{\alpha_1}(\theta (0,\alpha_2)) \quad \text{ for } k= 1\dots,K
 	\]
 \end{lemma}
 
 %In our proposed method, we use the dynamic programming method \citep{johnson2013dynamic} to solve the FLSA problem, which only runs in $O(K)$ time for the subproblem with $K$ dimensional variables.

\subsubsection*{B. Group lasso penalty $P_G$}

By definition, the update of $\Theta_{t+1}$ for group lasso penalty $P_G(\Theta)$ is:
\begin{equation*}
 \bm{\Theta}_{t+1}  =	\arg\min_{\Theta} \left\{ \frac{1}{2}||\bm{\Theta}-\bm{\Theta}_{t}+\eta_{t} \nabla f(\bm{\Theta}_{t}) ||^2 
  \eta_t \cdot \lambda_1 \sum_{k=1}^K \sum_{i\neq j}|\theta_{k,i,j}|+\eta_t \cdot \lambda_2 \sum_{i\not=j}(\sum_{k=1}^K \theta_{k,i,j}^2)^{\frac{1}{2}} \right\}
\end{equation*}

Similarly, let $A=\bm{\Theta}_{t}-\eta_{t} \nabla f(\bm{\Theta}_{t}) $; then, the problem becomes
\begin{equation*}
  \underset{\theta_{1,i,j},\cdots,\theta_{K,i,j}}{\text{argmin }}  	 \left\{ \frac{1}{2} \sum_{k=1}^K(\theta_{k,i,j}-a_{k,i,j}  )^2   + \eta_t \cdot \lambda_1 \cdot 1_{i\neq j} \sum_{k=1}^K |\theta_{k,i,j}|  +\eta_t \cdot \lambda_2 \cdot 1_{i\neq j} (\sum_{k=1}^K \theta_{k,i,j}^2)^{\frac{1}{2}}  
\right\} .
\end{equation*}

Obviously, when $i=j$, $\theta_{k,i,j}=a_{k,i,j}$. For $i\neq j$, it is a group lasso problem \citep{yuan2006model,friedman2010note}, and we can simply solve it by the following \cite{friedman2010note}:

\[
\theta_{k,i,j} = {\cal S}_{\eta_t \cdot \lambda_1} (a_{k,i,j}) \left ( 1-\frac{\eta_t \cdot \lambda_2}{\sqrt{\sum_{k=1}^{K} {\cal S}_{\eta_t \cdot \lambda_1} (a_{k,i,j})^2 }}  \right)_{+}.
\]

\subsubsection{Theoretical analysis}
In the following claim, we extend the theoretical analysis in Section 2 to the JGL problem and show that the optimal solution $\bm{\Theta}_*$ in the JGL problem can be bounded. 

As we can see in Lemma \ref{lm2}, the bounds of $\bm{\Theta}_t$ in Algorithm \ref{alg:algorithm0} are related to the bounds of $\bm{\Theta}_*$, which were already proven by \cite{banerjee2008model} and \cite{lu2009smooth}.
For multiple Gaussian graphical models, \cite{honorio2010multi} and \cite{hara2013learning} provided the bounds for the optimal solution $\bm{\Theta}_*^{(k)}$. However, though their target models were similar, they were different from those of the JGL model.
  To the best of our knowledge, no related research has provided the bounds of the optimal solution $\bm{\Theta}_*^{(k)}$ for the JGL problem.

In the following, we provide the lower and upper bounds for $\bm{\Theta}_*^{(k)}$, which are applied to both fused and group lasso-type penalties.

\begin{proposition} \label{propo2}
The optimal solution $\bm{\Theta}^*$ of the problem \eqref{eq6} satisfies
\[
  \frac{n_k}{p \lambda_c +n_k ||S^{(k)}||_2}
 \leq ||\bm{\Theta}_*^{(k)}||_2 \leq   \frac{Np}{\lambda_1} + \sum_{k=1}^K\sum_{i=1}^p (s_{ii}^{(k)})^{-1}
\]
where $\lambda_c :=\sqrt{K \lambda_1^2 + 2 K \lambda_1 \lambda_2  +\lambda_2^2}$.
\end{proposition}

We use the same G-ISTA strategy to handle the initial step size of $\eta_{t}$ for each iteration in Algorithm \ref{alg:algorithm1}. In the numerical experiment (section 4.2.3), the algorithm also shows a linear convergence rate. %

\subsection{Modified ISTA for JGL}

In the previous discussion, we mentioned that $f(\bm{\Theta})$ in the JGL problem is not globally Lipschitz gradient continuous. The ISTA may not be efficient enough for the JGL case because it includes the backtracking line search procedure for this case, which needs to evaluate the objective function and is inefficient when the evaluation is expensive. 

In this section, we modify Algorithm \ref{alg:algorithm1} to Algorithm \ref{alg:algorithm3} based on the strategy in Section \ref{cscm}, which takes advantage of the properties of the self-concordant function. The self-concordant function does not rely on the Lipschitz gradient assumption on the smooth part $f(\bm{\Theta})$ \citep{tran2015composite}, and we can eliminate the need for the backtracking line search. %

\begin{lemma}(\cite{boyd2004convex}) \label{sclm}
 Self-concordance is preserved by scaling and addition: if $f$ is a self-concordant function and a constant $a \leq 1$, then $af$ is self-concordant. If $f_1,f_2$ are self-concordant, then $f_1+f_2$ is self-concordant.
\end{lemma}

By Lemma \ref{sclm}, the function \eqref{ftheta} %\[f(\bm{\Theta})=\sum_{k=1}^{K} n_k \{  -\text{logdet} \bm{\Theta}^{(k)} +\text{trace}(\bm{S}^{(k)}\bm{\Theta}^{(k)}) \}\] 
is a self-concordant function. %Then, we discuss the convergence rate of the Algorithm \ref{alg:algorithm3} by using the properties of the self-concordant function $f(\Theta)$.
In Algorithm \ref{alg:algorithm3}, for the initial step size of $L_{t}^{-1}$ in each iteration, we use the same strategy as that of Algorithm \ref{alg:algorithm1}. 
Then, the mechanism in \cite{tran2015composite} is employed in Steps 3-5 of Algorithm \ref{alg:algorithm3}. 

\begin{algorithm}[H]
\caption{Modified ISTA algorithm}
\label{alg:algorithm3}
\textbf{Input}: $\bm{S}$, tolerance $\epsilon>0$, initial step size $L_{0}^{-1}$, initial iterate $\bm{\Theta}_0$.  \\
\textbf{For} $t=0,1,\cdots,  $ (until convergence) \textbf{do} 
\begin{algorithmic}[1] %[1] enables line numbers
\STATE Initialize $L_t$ 
\STATE Compute 
\begin{align}
\bm{d}_t:=\text{prox}_{L_{t}^{-1} g} (\bm{\Theta}_{t}- L_{t} ^{-1}\nabla f(\bm{\Theta}_{t}))-\bm{\Theta}_{t}
\end{align}
\STATE Compute 
    $\beta_t:=L_t||vec(\bm{d}_t)||_{2}^2 $
and   
    $ \lambda_t:= \sum_{k=1}^{K} n_k || (\bm{\Theta}_t^{(k)})^{-1}\bm{d}_t^{(k)}||_{2}.$
\STATE Determine step size $\alpha_t:=   \frac{\beta_t}{\lambda_t(\lambda_t+\beta_t)}$
\STATE If $\alpha_t>1$, then set $L_t := \frac{L_t}{2}$ and go back to Step 2
\STATE Update $\bm{\Theta}_{t+1}:=\bm{\Theta}_{t}+\alpha_t \bm{d}_t$
\end{algorithmic}

\textbf{Output}: $\epsilon$-optimal solution to problem (1), $\bm{\Theta}^*=\bm{\Theta}_{t+1 }$
\end{algorithm}

However, there is no backtracking procedure in this algorithm that guarantees the positive definiteness of $\bm{\Theta}_t$, as in \eqref{condi} of Algorithm \ref{alg:algorithm1}. Hence, in Algorithm \ref{alg:algorithm3}, we need to illustrate how to ensure the positive definiteness of the $\bm{\Theta}_{t+1}$ in the iterations.

\begin{lemma}(\cite{nemirovski2004interior}, Theorem 2.1.1) \label{epp}
	Let $f$ be a self-concordant function and let $x \in dom f$. Additionally, if
	\[
	W_r(\bm{x}) = \{ \bm{y}| (\bm{y}-\bm{x}) \nabla^2 f(\bm{x}) (\bm{y}-\bm{x}) \leq 1 \}
	\]
	then $W_r(\bm{x}) \subset dom f$.
\end{lemma}

In Algorithm \ref{alg:algorithm3}, because we know $\alpha_t:=\frac{\beta_t}{\lambda_t (\lambda_t +\beta_t)}$ and $\alpha_t \leq 1$ by Step 3 and 5, therefore $\alpha_t^2 \lambda_t <1$, i.e.:
\[
 \alpha_t^2 \lambda_t = \alpha_t^2  \sum_{k=1}^{K} n_k || (\bm{\Theta}_t^{(k)} )^{-1}\bm{d}_t^{(k)}||_2 <1
\]
\[
(\bm{\Theta}_{t+1}- \bm{\Theta}_t) \nabla^2 f(\bm{\Theta}_t) (\bm{\Theta}_{t+1}- \bm{\Theta}_t) <1
\]
Hence, by Lemma \ref{epp}, the $\bm{\Theta}_{t+1}$ stays in the domain and keeps the positive definiteness.

By Lemma 12 in \cite{tran2015composite}, it is known that the function value is always decreasing in Algorithm \ref{alg:algorithm3}:
\begin{equation} \label{desc}
  F(\bm{\Theta}_{t+1}) \leq F(\bm{\Theta}_t) \qquad \text{for } t = 0,1,\dots,
\end{equation}

Then, based on condition \eqref{desc}, we provide the explicit bounds of $\{ \bm{\Theta}_t \}_{t=0,1\dots,}$ in Algorithm \ref{alg:algorithm3} for the JGL problem with the help of Lemma \ref{lm2} and Proposition \ref{propo2}. For the proof, see Appendix.

\begin{proposition}
Sequence $\{\bm{\Theta}_t\}_{t=0,1,\cdots,}$ generated by Algorithm \ref{alg:algorithm3} can be bounded: 
	\[
	m \leq ||\bm{\Theta}_t||_2 \leq M 
	\]
	where $M:= ||\bm{\Theta}_{0}||_F + 2||\bm{\Theta}_*||_F,m:=e^{-\frac{C_1}{n_m}} M^{(1-Kp)} $, $n_m = \max_k{n_k}$,and constant $C_1:=F(\bm{\Theta}_0)$.
\end{proposition}

\section{Experiments}

In this section, we evaluate the performance of the proposed method on both synthetic and real datasets, and we compare the following algorithms:
\begin{itemize}
    \item ADMM: the general ADMM method proposed by \cite{danaher2014joint}
    \item FMGL: the proximal Newton-type method proposed by \cite{yang2015fused}.
    \item ISTA: the proposed method in Algorithm \ref{alg:algorithm1}.
    \item M-ISTA: the proposed method in Algorithm \ref{alg:algorithm3}.
\end{itemize}

We perform all the tests in R Studio on a Macbook Air with 1.6 GHz Intel Core i5 and 8 GB memory. The wall times are recorded as the run times for the four algorithms.

\subsection{Stopping criteria}

In the experiments, we consider two stopping criteria for the algorithms.

1. Relative error:
\[
\frac{ \sum_{k=1}^K ||\bm{\Theta}_{t+1}^{(k)}-\bm{\Theta}_{t}^{(k)}||_{F}  } {\max \{ \sum_{k=1}^K ||\bm{\Theta}_{t}^{(k)}||_{F},1 \}  } \leq  \epsilon
\]

2. Objective error:
\[
 ||F(\bm{\Theta}_{t})-F(\bm{\Theta}_{*})||_{F}   \leq  \epsilon
\]

We use the objective error for convergence rate analysis and the relative error for the time comparison. Because, in our proposed method, we do not obtain dual solutions, we do not employ the duality gap stopping criterion in Section \ref{istagl}.

%$||\Theta_{i+1}-\Theta_{i} ||_F \leq  10^{-8} \times \max \{ ||\Theta_{i+1}||_F,1 \} $
\subsubsection{Model selection}
%bayesian information criterion(BIC) or through crossvalidation
The JGL model is affected by regularized parameters $\lambda_1$ and $\lambda_2$.
For selecting the parameters, we use the D-fold cross-validation method. First, the dataset is randomly split into $D$ segments of equal size, a single subset (test data), estimated by the other $D-1$ subsets (training data) and change the subset for the test to repeat $D$ times so that each subset is used.

Let $S^{(k)}_d$ be the sample covariance matrix of d-th ($d=1,\ldots,D$) segment for class $k=1,\dots,K$.
We estimate the inverse covariance matrix by the remaining $D-1$ subsets $\hat{ \Theta}^{(k)}_{\lambda, -d}$, and choose $\lambda_1$ and $\lambda_2$ that minimize the average predictive negative log-likelihood as follows:
\[ 
CV( \lambda_1,\lambda_2 ) =  \sum_{d=1}^D \sum_{k=1}^K \left\{ n_k \text{trace}(\bm{S}^{(k)}_d \hat{\bm{\Theta}}^{(k)}_{\lambda, -d} ) - \text{logdet} \hat{ \bm{\Theta}}^{(k)}_{\lambda, -d} \right\} .
\]

\subsection{Synthetic data}

%We control the number of non-zero entries (nnz) in each $\Theta^{(k)}$ to be about 10p, so that  the total number of nonzeros in the K precision matrices is 10Kp. 

 We follow the data generation mechanism described in \cite{li2006gradient} with some modifications for the JGL model. We put the details in Appendix.

The performance of the proposed methods was assessed on synthetic data in terms of the number of iterations, the execution time, the squared error, and the receiver operating characteristic (ROC) curve.

\subsubsection{Comparison experiments}
\begin{figure}[H]
  \begin{center}
  \includegraphics[width=3.7in,height=3.7in]{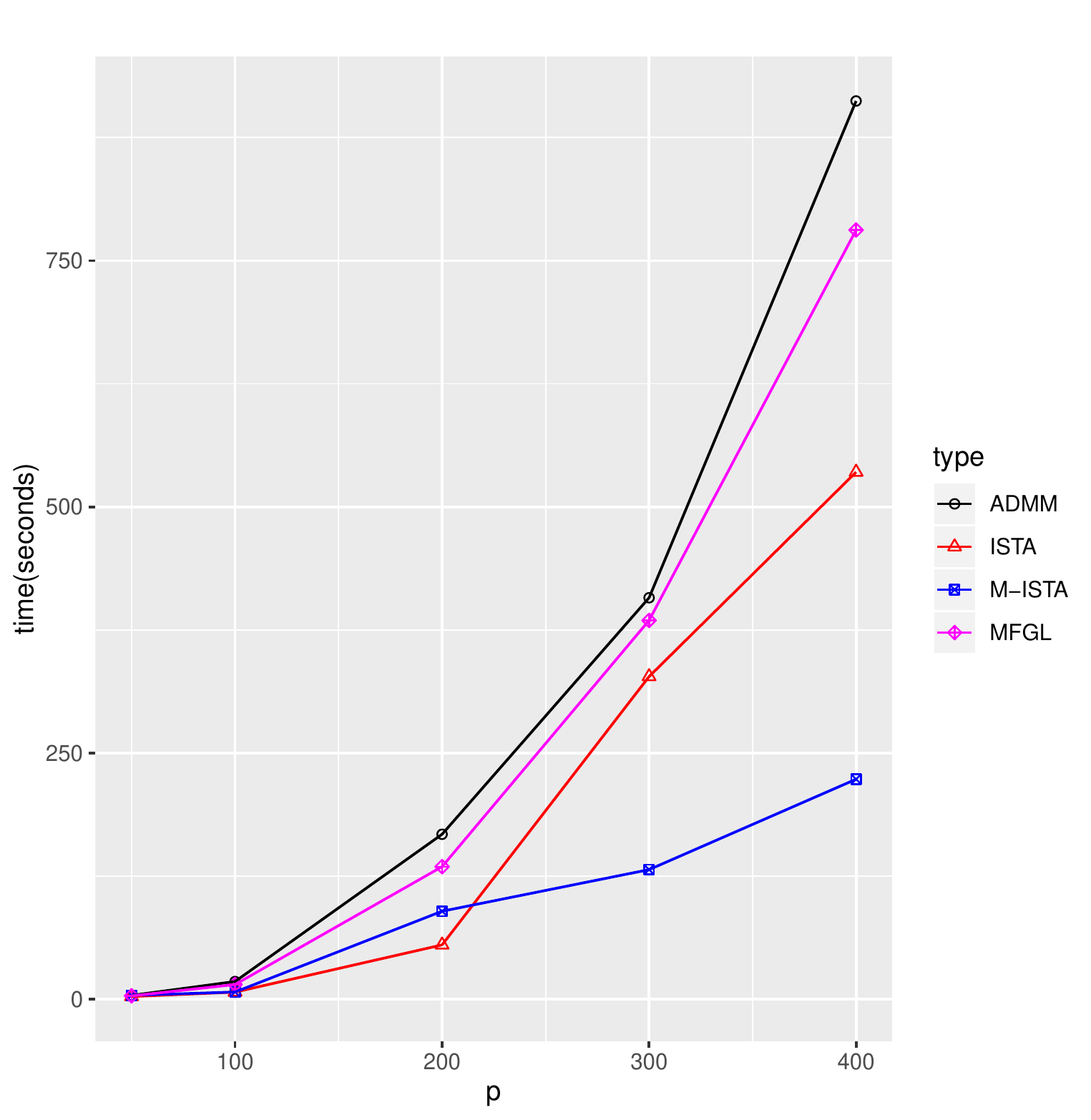}
    \caption{Plot of time comparison under different $p$. Setting $\lambda_1=0.1,\lambda_2=0.05, K=2$ and $N=200$.}\label{compp}
  \end{center}
\end{figure}

We vary $p, N, K \text{ and } \lambda_1$ to compare the execution time of our proposed methods with that of the existing methods. For a fair comparison in the experiments, we consider only the fused penalty in our proposed method because the FMGL algorithm applies only to the fused penalty.
First, we compare the performance among different algorithms under various dimensions $p$, which are shown in Figure \ref{compp}.

Figure \ref{compp} shows that the execution time of the FMGL and ADMM increases	rapidly as $p$ increases. In particular, we observe that the M-ISTA significantly outperforms when $p$ exceeds 200. The ISTA shows better performance than the three methods when $p$ is less than 200, but it requires more time as $p$ grows, compared to the M-ISTA. It is reasonable to consider that evaluating the objective function in the backtracking line search at every iteration increases the computational burden, especially when $p$ increases, which means that the M-ISTA is a good choice for these cases. Furthermore, the ISTA can be a good candidate when the evaluation is inexpensive.

\begin{table}[H]
\caption{Computational time under different settings \label{tab02}}
\begin{tabular}{|l|l|l|l|l|l|l|l|l|l|}
\hline
\multicolumn{6}{|l|}{Parameters setting} &
  \multicolumn{4}{l|}{Computational time} \\ \hline
$p$ &
  $K$ &
  $N$ &
  1 &
  2 &
  precision $\epsilon$ &
  ADMM &
  FMGL &
  ISTA &
  M-ISTA \\ \hline
\multirow{2}{*}{20} &
  2 &
  \multirow{2}{*}{50} &
  0.1 &
  0.05 &
  \multirow{2}{*}{0.00001} &
  10.506 secs &
  \textbf{1.158 secs} &
  2.174 secs &
  1.742 secs \\ \cline{2-2} \cline{4-5} \cline{7-10} 
 &
  5 &
   &
  1 &
  0.5 &
   &
  1.123mins &
  10.556 secs &
  4.216 secs &
  \textbf{2.874 secs} \\ \hline
\multirow{2}{*}{30} &
  2 &
  \multirow{2}{*}{120} &
  0.1 &
  0.05 &
  \multirow{2}{*}{0.0001} &
  36.592 secs &
  13.225 secs &
  \textbf{12.675 secs} &
  12.700 secs \\ \cline{2-2} \cline{4-5} \cline{7-10} 
 &
  3 &
   &
  0.1 &
  0.05 &
   &
  3.779 mins &
  2.424 mins &
  \textbf{58.208 secs} &
  1.481 mins \\ \hline
\multirow{3}{*}{50} &
  \multirow{3}{*}{2} &
  \multirow{3}{*}{600} &
  0.02 &
  \multirow{3}{*}{0.005} &
  \multirow{3}{*}{0.0001} &
  6.427 secs &
  10.228 secs &
  7.213 secs &
  \textbf{4.625 secs} \\ \cline{4-4} \cline{7-10} 
 &
   &
   &
  0.03 &
   &
   &
  6.240 secs &
  8.925 secs &
  6.645 secs &
  \textbf{4.023 secs} \\ \cline{4-4} \cline{7-10} 
 &
   &
   &
  0.04 &
   &
   &
  7.025 secs &
  9.381 secs &
  6.144 secs &
  \textbf{3.993 secs} \\ \hline
\multirow{3}{*}{200} &
  \multirow{3}{*}{2} &
  \multirow{3}{*}{400} &
  0.09 &
  \multirow{3}{*}{0.05} &
  \multirow{3}{*}{0.0001} &
  4.050 mins &
  1.874 mins &
  2.289 mins &
  \textbf{35.038 secs} \\ \cline{4-4} \cline{7-10} 
 &
   &
   &
  0.1 &
   &
   &
  4.569 mins &
  1.137 mins &
  1.340 mins &
  \textbf{24.852 secs} \\ \cline{4-4} \cline{7-10} 
 &
   &
   &
  0.12 &
   &
   &
  3.848 mins &
  1.881 mins &
  1.443 mins &
  \textbf{18.367 secs} \\ \hline
\end{tabular}
\end{table}

Table \ref{tab02} summarizes the performance of the four algorithms under different parameter settings for achieving a given precision $\epsilon$ of the relative error.
The results presented in Table \ref{tab02} reveals that when we increase the number of classes $K$, all the algorithms spend more time than usual. Moreover, the execution time of ADMM becomes huge among them.
When we vary the $\lambda_1$, they become more efficient as the value larger.
For most instances, the M-ISTA and ISTA outperform the existing ones, such as ADMM and FMGL. As for the exceptional cases ($p=20, k=2, N=50, \lambda_1 = 0.1$ and $\lambda_2 = 0.05$), the M-ISTA is still comparable with the FMGL and faster than ADMM.

\subsubsection{Algorithm Assessment}
%A number of metrics are used to assess performance, including receiver operating characteristic curves, average entropy loss, average Frobenius loss, average false positive and avergae false negative rates and the average rate of misidentified common zeros among the categories,

We assessed our proposed method by drawing a receiver operating characteristic (ROC) curve, which displays the number of true positive edges (i.e., TP edges) selected against the number of false edges  (i.e., FP edges) selected. We say that an edge $(i,j)$ in the $k$-th class is selected in estimate $\hat{\bm{\Theta}}^{(k)}$ if element $\hat{\theta}^{(k)}_{ij}\neq 0$ and the edges are true positive edges selected if the precision matrix $\theta^{(k)}_{ij}\neq 0$ and false positive edges selected if the precision matrix $\theta^{(k)}_{ij}=0$, where the two quantities are defined by
\[
TP= \sum_{k=1}^K \sum_{i<j} 1 (\theta_{ij}^{(k)}=0)  \cdot 1 (\hat{\theta}_{ij}^{(k)}=0)\ ,
\]
\[
FP= \sum_{k=1}^K \sum_{i<j} 1 (\theta_{ij}^{(k)} \neq 0)  \cdot 1 (\hat{\theta}_{ij}^{(k)}=0)\ ,
\]
and $1(\cdot)$ is the indicator function.

To confirm the validity of the proposed method, we compare the ROC figures of the fused penalty and group penalty. We fix the parameters $\lambda_2$ for each curve and change the $\lambda_1$ value to obtain various numbers of selected edges because the sparsity penalty parameter $\lambda_1$ can control the number of selected total edges.

\begin{figure}[htbp]
\centering
\subcaptionbox{The fused penalty.\label{comproc}}[.48\linewidth]
{
	\includegraphics[scale = 0.7 ]{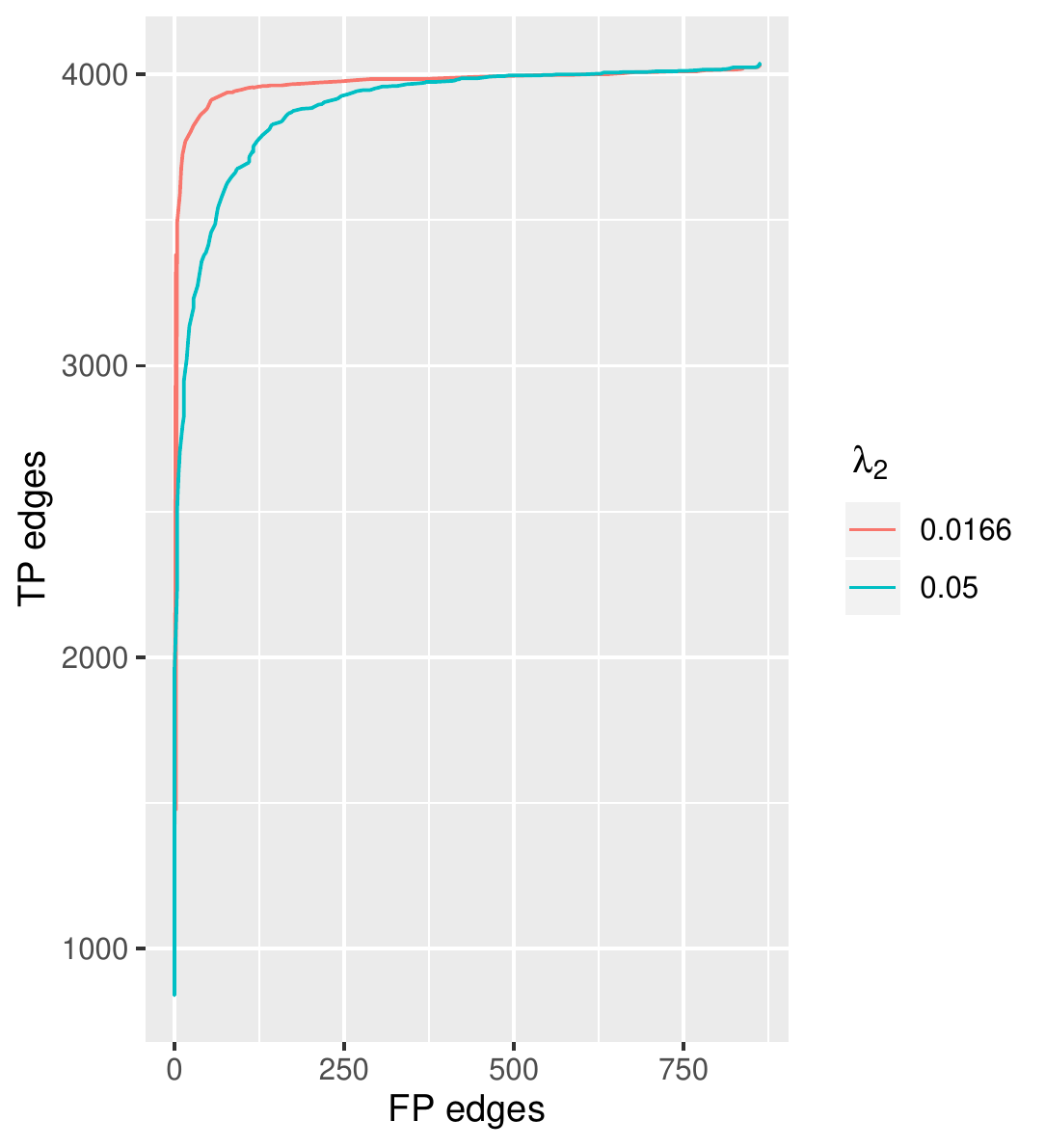}
}
\subcaptionbox{The group penalty.\label{comprocg}}[.48\linewidth]
{
	\includegraphics[scale = 0.7 ]{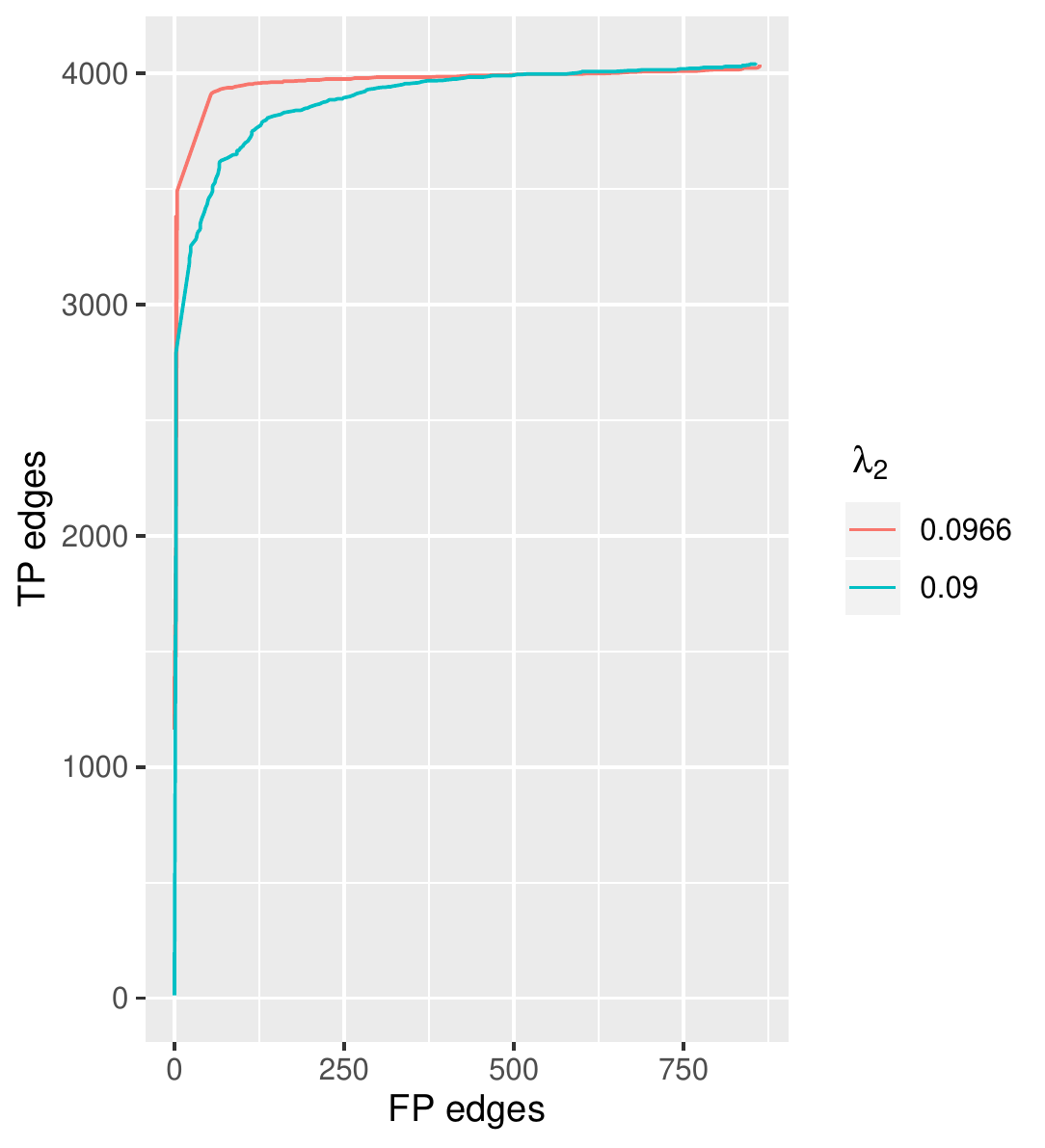}
}
\caption{Plot of true positive edges vs. false positive edges selected. Setting p=50, K=2.}
\end{figure}

%Receiver operating characteristic curves. The horizontal and vertical axes ineach panel are false positive rate and TP rate (sensitivity), respectivityly.
First, we show the ROC curves for fused and group lasso penalties in Figure \ref{comproc} and Figure \ref{comprocg}, respectively.

From the figures, we observe that both penalties show highly accurate predictions for the edge selections. The result of $\lambda_2=0.0166$ in the fused penalty case is better than that in $\lambda_2 = 0.05$. Additionally, the result of $\lambda_2=0.0966$ in the fused penalty case is better than that in $\lambda_2 = 0.09$, which means that if we select the tuning parameters properly, then we can obtain precise results while simultaneously meeting our different model demands. 

\begin{figure}[htbp]
\centering
\subcaptionbox{The fused penalty.\label{compmse}}[.48\linewidth]
{
	\includegraphics[scale = 0.7 ]{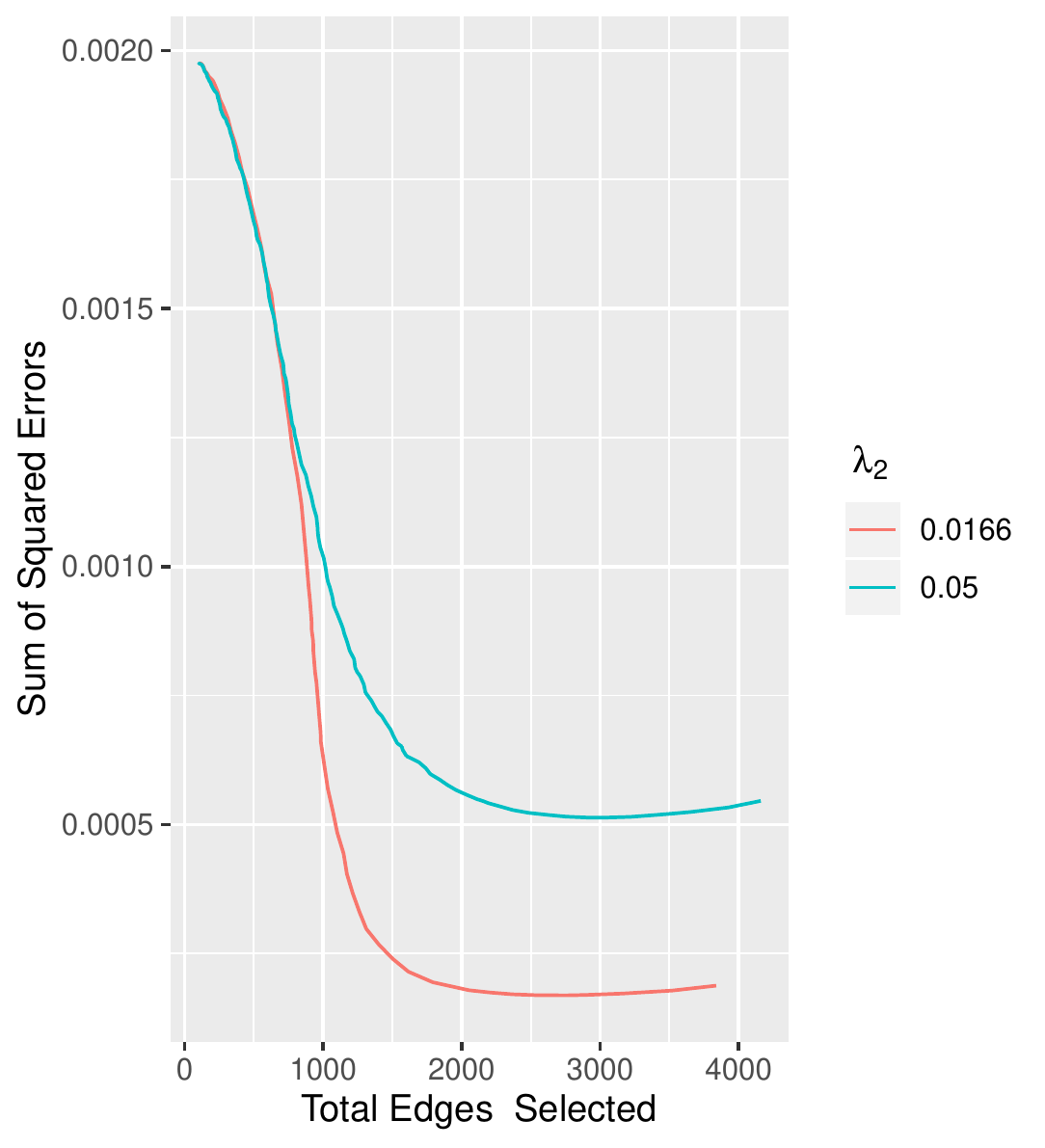}
}
\subcaptionbox{The group penalty.\label{compmseg}}[.48\linewidth]
{
	\includegraphics[scale = 0.7 ]{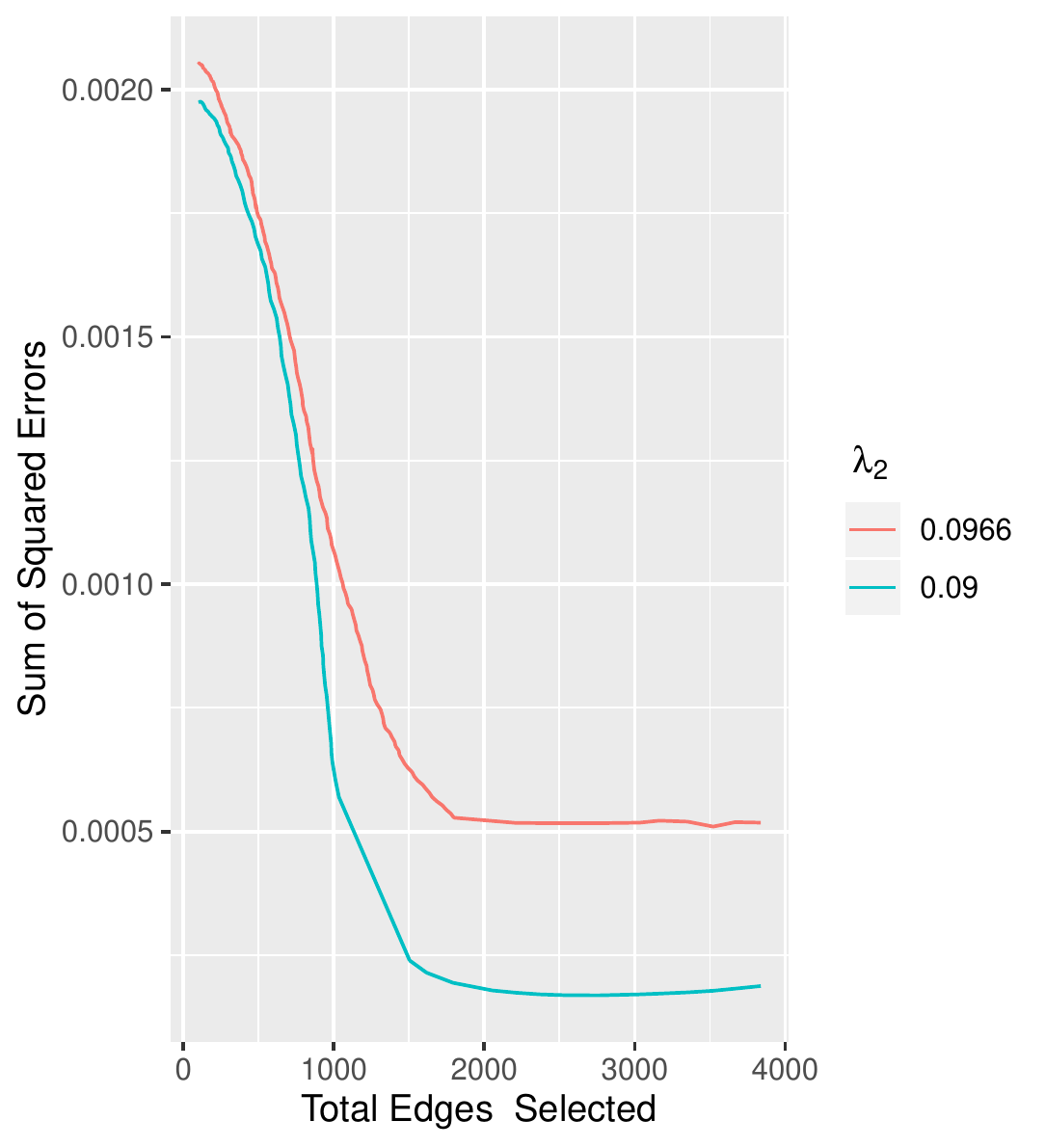}
}
\caption{Plot of the mean squared errors vs. total edges selected. Setting p=50, K=2.}
\end{figure}

Then, Figure \ref{compmse} and Figure \ref{compmseg} display the mean squared errors (MSE) between the estimated values and true values. 

\[
MSE = \frac{2}{Kp(p-1)} \sum_{k=1}^K \sum_{i<j} ( \hat{\theta}_{ij}^{(k)} -{\theta}_{ij}^{(k)} )^2
\]
where $\hat{\theta}_{ij}^{(k)}$ is the value estimated by the proposed method, and ${\theta}_{ij}^{(k)}$ is the true precision matrix value we use in the data generation.

The figures illustrate that when the total number of edges selected is increasing, the errors decrease and finally achieve relatively low values.

Overall, the proposed method shows competitive efficiency not only in computational time but also in  accuracy.

\subsubsection{Convergence rate}
This section shows the convergence rate of the ISTA for the JGL problem in practice, with $\lambda_1=0.1, 0.09$ and $0.08$. We recorded the number of iterations to achieve the different tolerance of $F(\bm{\Theta}_t)-F(\bm{\Theta}_*)$ and ran it on a synthetic dataset, with $p=200, K=2,\lambda_2 = 0.05$ and $N=400$. The figure reveals that as $\lambda_1$ decreases, more iterations are needed to converge to the specified tolerance. Moreover, the figure shows the linear convergence rate of the proposed ISTA method.

\begin{figure}[H]
  \begin{center}
  \includegraphics[width=4in,height=5in]{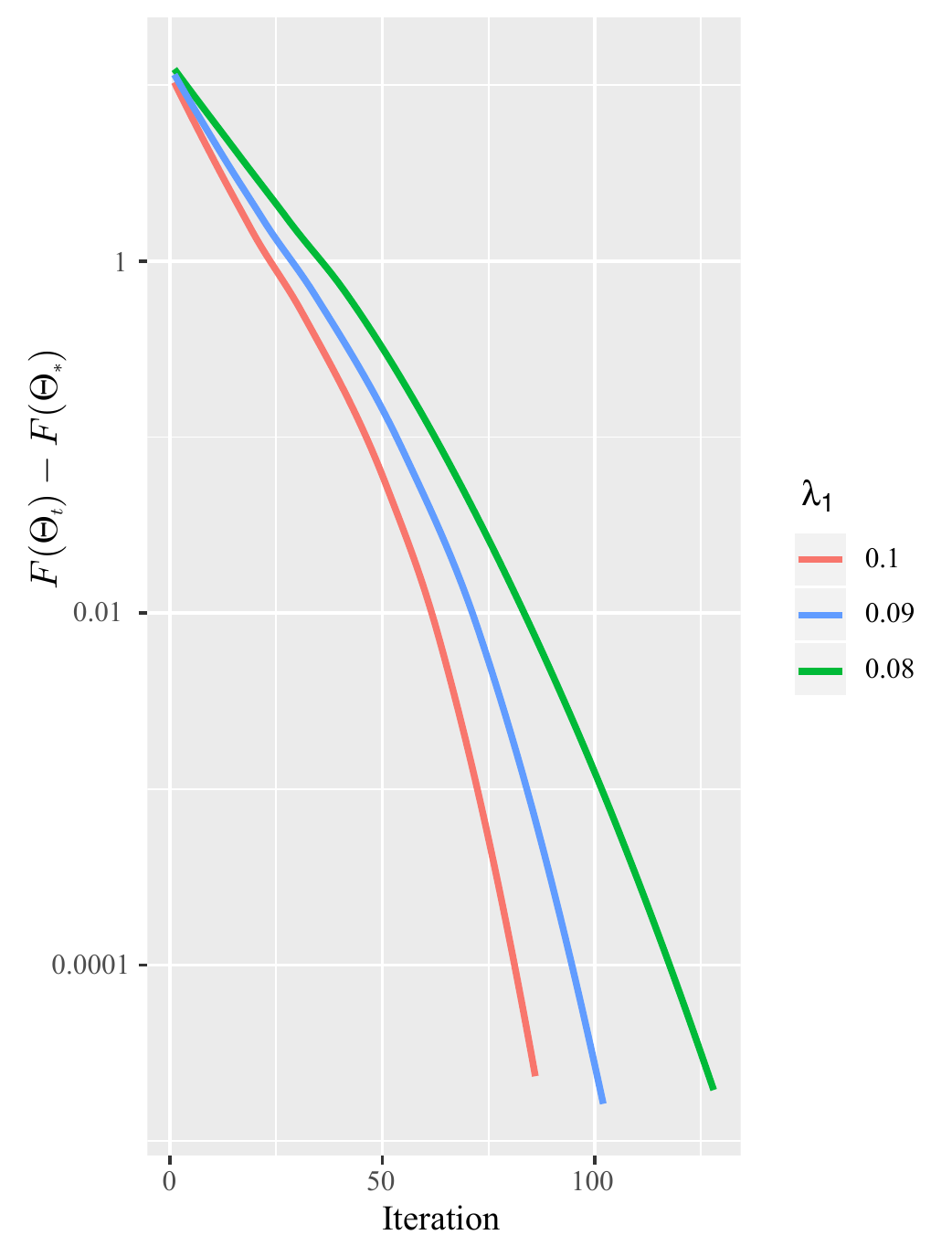}
    \caption{Plot of $\log ( F(\bm{\Theta}_t) - F(\bm{\Theta}_*) )$ vs. the number of iterations with different $\lambda_1$ values. Setting $p=200, N= 400,K=2$ and $\lambda_2=0.05$.}\label{convera1}
  \end{center}
\end{figure}

\subsection{Real data}
We use a breast cancer dataset for the experiment to show the common structure and jointly estimate common links across graphs. There are 250 samples and 1,000 types of genes in the dataset, with 192 control samples and 58 case samples. Additionally, we set $\lambda_1= 0.3,\lambda_2= 0.08$ for suitable visualization. %the data source?
\begin{figure}[H]
\begin{center}
	\includegraphics[width=\textwidth]{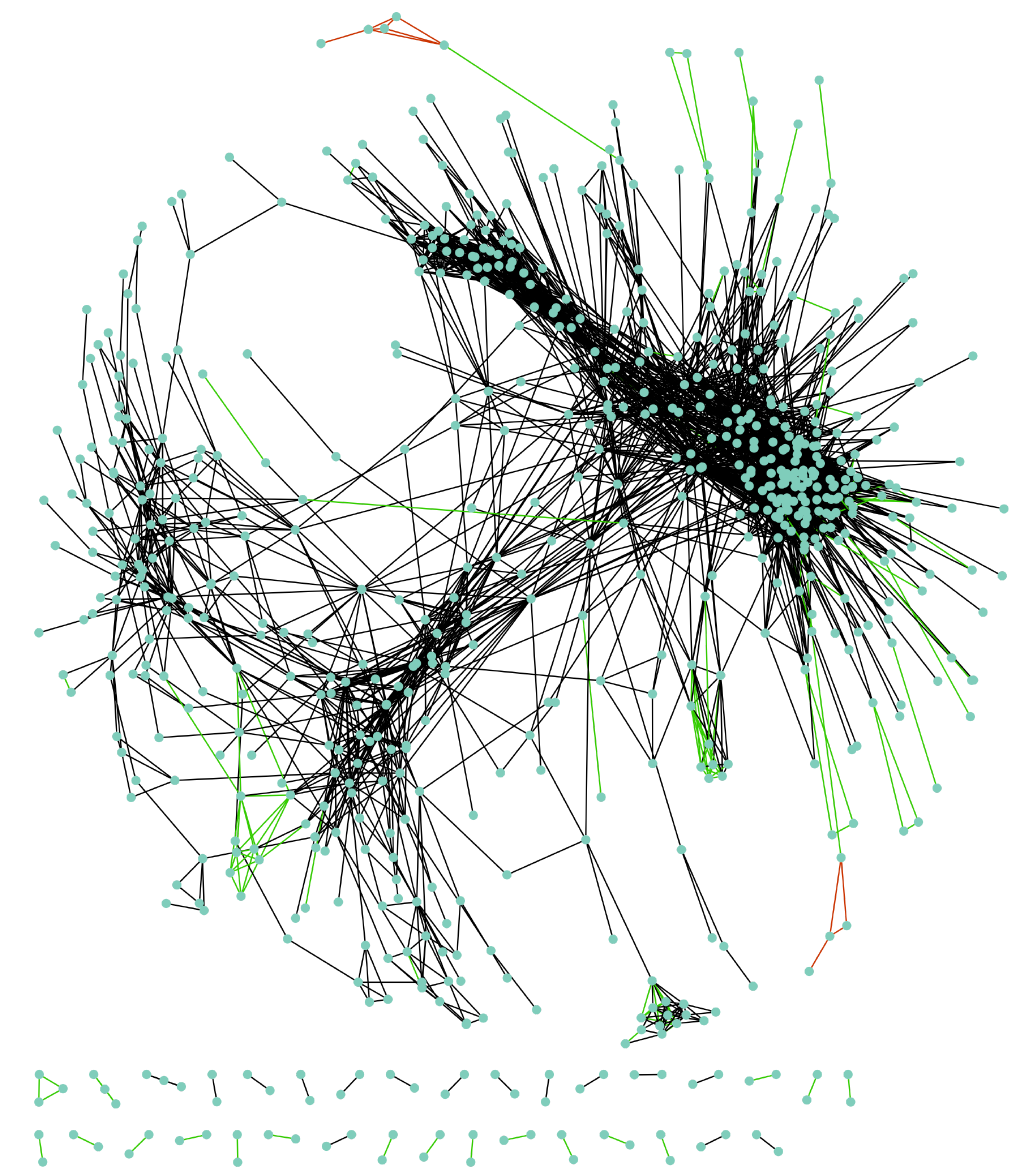}
\end{center}
	\caption{Graph of the gene expression of breast cancer} \label{bcf}
\end{figure}

We use Cytoscape \citep{shannon2003cytoscape} to visualize the results. Figure \ref{bcf} shows the graph of the breast cancer gene expression identified by the JGL model with the proposed methods. As we can see in the figure, each node represents a gene, and the edges demonstrate the relationship between genes.

We use different colors to show various structures. The green edges are a common structure between the two classes, the red edges are the specific structures for the control classes, and the black edges are for the case samples.

\section{Conclusions}
We propose two efficient proximal gradient descent procedures with and without the backtracking line search option for the joint graphical lasso. The first does not require extra variables, unlike ADMM, which needs manual tuning of the Lagrangian parameters and dual variables. Moreover, we reduce it to subproblems that can be solved both efficiently and precisely. The second does not require backtracking line search, which significantly reduces the computation time needed to evaluate objective functions.

From the theoretical perspective, we reach the linear convergence rate for the ISTA. Furthermore, we derive the lower and upper bounds of the solution to the JGL problem and the iterations in the algorithms, guaranteeing that the iterations stay in the constrained domain. Numerically, the methods are demonstrated on both simulated and real datasets to illustrate their robust and efficient performance over state-of-the-art algorithms.

For the further computational improvement, the most expensive step in the algorithms is to calculate the inversion of matrices required by the gradient of $f(\bm{\Theta})$. Both algorithms have A complexity of $O(Kp^3)$ per iteration. Moreover, we can solve the matrix inversion problem with more efficient algorithms with lower complexity. In addition, we can also use the faster computation procedure in \cite{danaher2014joint} to decompose the optimization problem for the proposed methods and regard it as preprocessing. Overall, the proposed methods are highly efficient for the joint graphical lasso problem.

\newpage

\bibliographystyle{plainnat}
\bibliography{Bibliography-MM-MC}

\end{document}